\documentclass{article}
\usepackage[cp1251]{inputenc}
\usepackage[english]{babel}
\usepackage{amsmath}
\usepackage{amsfonts,amssymb,mathrsfs,amscd}
\usepackage{color}

\newtheorem{theorem}{Theorem}
\newtheorem{lemma}{Lemma}

\begin{document}


{\Large On the convergence of formal power-log series solutions of an algebraic ODE}
\bigskip

\centerline{R.\,R.\,Gontsov,\;I.\,V.\,Goryuchkina}

\begin{abstract}
We propose a sufficient condition of the convergence of a power-log series that formally satisfies an algebraic ordinary differential equation (ODE) of arbitrary order. A general form and properties of the functional coefficients of such a series are established.
\end{abstract}


\section{Introduction}\label{s1}

We study the question of the convergence of formal series in the powers of the independent variable $x\in\mathbb C$ with  coefficients rational in $\ln x$,
\begin{equation}\label{series}
\varphi=\sum_{k=0}^{\infty}R_k(\ln x)\,x^{k}, \qquad R_k\in{\mathbb C}(t),
\end{equation}
that satisfy an {\it algebraic} ODE of order $n$,
\begin{equation}\label{ODE}
F(x,y,\delta y,{\dots},\delta^n y)=0,
\end{equation}
where $F=F(x,y_0,y_1,\ldots,y_n)$ is a polynomial of $n+2$ variables and $\delta$ is the differentiation $x(d/dx)$.

Such series will be referred to as {\it power-log series}. One may regard them as the generalization of
{\it Dulac series} (whose coefficients $R_k$ are polynomials) and expect that a sufficient condition of the convergence of a formal power-log series solution is similar to that of a formal Dulac series solution from the paper \cite{GG}. Here we use again a universal majorant method but the proof in the present article also contains quite a few new moments  compared to \cite{GG}: this is in part due to the fact that the coefficients of a formal Dulac series solution of (\ref{ODE}) that satisfies the sufficient condition of convergence from \cite{GG}, starting with some number $k$ are solutions of nonhomogeneous linear ODEs of order $n$ with {\it constant} coefficients, whereas the coefficients of power-log series considered here are solutions of nonhomogeneous linear ODEs of order $n$ with {\it rational} coefficients. In particular, we will show that the poles of all the coefficients $R_k$ of the series (\ref{series}) that formally satisfies (\ref{ODE}), belong to a {\it common} finite set and that their orders have a linear growth with respect to $k$.

Power-log series (\ref{series}) form a subclass in a more general class of formal series in the powers of the independent variable $x$ with coefficients in the form of {\it formal} Laurent series in $1/\ln x$ with a finite main part. Elements of such a general class are sometimes referred to as {\it power-log transseries} \cite{Ed}. They are rather often found among formal solutions of nonlinear algebraic ODEs (for example, among formal solutions of the sixth Painlev\'e equation, see \cite{BG}, \cite{Guz}). Nevertheless, a formal, asymptotic and analytic theory of such series are at a starting point of development. At the moment, there are pretty of open questions. One of them is the following: {\it under which conditions on the equation, the coefficients of its formal power-log transseries solution are rational functions of logarithm?} In our understanding, this question is more complicated than that we consider in the paper. There is no answer to this question yet, there are only some special  cases where one can prove the rationality (in $\ln x$) of the coefficients of power-log transseries. For example, using the relation of the sixth Painlev\'e equation with the Schlesinger equation of isomonodromic deformations of Fuchsian systems, S.\,Shimomura \cite{Sh} has obtained such a rationality for some of formal power-log transseries solutions of the sixth Painlev\'e equation and proved the convergence of the formal solutions themselves (see also Example 1 at the end of our article).
So, stress again that in this article we assume that the equation (\ref{ODE}) {\it possesses} a formal power-log series solution (\ref{series}). Then we propose the following general condition of its convergence.

\begin{theorem}\label{th1}
Let the series $(\ref{series})$ formally satisfy the equation $(\ref{ODE})$:
$$
F(x,\Phi):=F(x,\varphi,\delta\varphi,\dots,\delta^n \varphi)=0,
$$
and let for each $j=0,\dots,n$ one have
\begin{equation}\label{partderiv}
\frac{\partial F}{\partial y_j}(x,\Phi)=a_j(\ln x) x^m+a_j^1(\ln x) x^{m+1}+\dots, \quad a_j, a_j^1,\ldots\in{\mathbb C}(t),
\end{equation}
where the number $m\in{\mathbb Z}_+$ is the same for all $j$, and $a_n\not\equiv0$.

Then for any open sector $S\subset\mathbb C$ of sufficiently small radius with vertex at the origin and of opening less than $2\pi$, the series $\varphi$ converges uniformly in $S$.
\end{theorem}

Theorem \ref{th1} is finally proved in the last Section \ref{s5} which is preceded by Sections \ref{s2},\ref{s3},\ref{s4} containing necessary auxiliary lemmas. After the proof of the theorem we present and discuss some examples which illustrate its possible applications.

\section{Reducing the ODE to a special form}\label{s2}

We begin the proof of Theorem 1 with a standard lemma on the reduction of the initial equation (\ref{ODE}) to a special form. The method of the proof of this lemma is adapted to formal solutions we study here but in essence it is similar to the method of the proof of the corresponding reduction lemma \cite{Ma} for an ODE possessing a formal power series solution (see also \cite[Lemma 5.1]{SS}).

\begin{lemma}\label{l1}
Under assumptions Theorem \ref{th1}, for any $\ell>m$ a transformation
$$
y=\sum\limits_{k=0}^{\ell}R_k(\ln x)x^k+x^\ell\, u
$$
reduces the initial equation $(\ref{ODE})$ to an equation of the form
\begin{equation}\label{rODE}
L(\delta)u=x\,M(x, \ln x, u,\delta u,\dots,\delta^n u),
\end{equation}
where
$$
L(\delta)=\sum\limits_{j=0}^n a_j(\ln x)(\delta+\ell)^j,
$$
and $M$ is a rational function with respect to $\ln x$ and polynomial with respect to the other $n+2$ variables.
\end{lemma}

{\bf Proof.} For any integer $\ell\geqslant0$, the formal series $\varphi$ can be presented in the form
$$
\varphi=\sum\limits_{k=0}^{\ell}R_k(\ln x)x^k+x^{\ell}\sum\limits_{k=1}^{\infty}R_{\ell+k}(\ln x)x^k=:\varphi_{\ell}+x^{\ell}\psi,
$$
whence
$$
\Phi=(\varphi,\delta\varphi,\ldots,\delta^n\varphi)=\Phi_{\ell}+x^{\ell}\Psi,
$$
where $\Phi_{\ell}=(\varphi_{\ell},\delta\varphi_{\ell},\ldots,\delta^n\varphi_{\ell})$ and $\Psi=(\psi,(\delta+\ell)\psi,\ldots,(\delta+\ell)^n\psi)$. Applying the Taylor formula one obtains
\begin{eqnarray}\label{Taylor}
0=F(x,\Phi_{\ell}+x^{\ell}\Psi)=F(x,\Phi_{\ell})+x^{\ell}\sum_{j=0}^n\frac{\partial F}{\partial y_j}(x,\Phi_{\ell})\psi_j+ \nonumber \\
+\frac{x^{2\ell}}2\sum_{i,j=0}^n\frac{\partial^2 F}{\partial y_i\partial y_j}(x,\Phi_{\ell})\psi_i\psi_j+\ldots,
\end{eqnarray}
where $\psi_j=(\delta+\ell)^j\psi$.
\medskip

Define the {\it valuation} of an arbitrary series
$$
\hat\varphi=\sum_{k=0}^{\infty}\widehat R_k(\ln x)x^k, \qquad \widehat R_k\in{\mathbb C}(t),
$$
as ${\rm val}(\hat\varphi):=\min\{k\mid \widehat R_k\not\equiv0\}$.
\medskip

Again, from the Taylor formula one obtains
$$
\frac{\partial F}{\partial y_j}(x,\Phi)-\frac{\partial F}{\partial y_j}(x,\Phi_{\ell})=x^{\ell}\sum_{i=0}^n
\frac{\partial^2 F}{\partial y_i\partial y_j}(x,\Phi_{\ell})\psi_i+\ldots,
$$
moreover ${\rm val}(\psi_i)\geqslant1$ for all $i$. From this and from (\ref{partderiv}) it follows that
$$
\frac{\partial F}{\partial y_j}(x,\Phi_{\ell})=a_j(\ln x) x^m+\tilde a_j^{(1)}(\ln x) x^{m+1}+\dots
$$
for each $j=0,1,\ldots,n$, that is, the leading coefficient $a_j(\ln x)$ is preserved when one substitutes the finite sum  $\Phi_{\ell}$ instead of $\Phi$ into $F'_{\partial y_j}$. From the relation (\ref{Taylor}) it follows that
$$
{\rm val}(F(x,\Phi_{\ell}))\geqslant m+\ell+1.
$$
Dividing (\ref{Taylor}) by $x^{m+\ell}$ we finally obtain an equation of the required form (\ref{rODE}), which possesses the formal solution
$$
\psi=\sum\limits_{k=1}^{\infty}R_{\ell+k}(\ln x)x^k.
$$
The lemma is proved.
\medskip

{\bf Remark 1.} Under assumptions of Theorem 1, without lost of generality one may put
$a_n(t)\equiv1$ in $L$. As for the right hand side of (\ref{rODE}), its function $M$ is a finite sum of the form
\begin{equation}\label{sumM}
M(x, \ln x, u,\delta u,\dots,\delta^n u)=\sum_{{\bf q},\mu}
x^{\mu}b_{{\bf q},\mu}(\ln x)u^{q_0}(\delta u)^{q_1}\ldots(\delta^nu)^{q_n},
\end{equation}
$\mu, q_j\in{\mathbb Z}_+$, $b_{{\bf q},\mu}\in{\mathbb C}(t)$.
\medskip

\begin{lemma}\label{l2}
The rational functions $R_{\ell+k}(t)$ determining the coefficients of the formal solution $\psi$ of $(\ref{rODE})$ can only have poles at points of a common finite set $\Delta=\{t_1,\ldots,t_s\}$ and maybe at infinity. For the order ${\rm ord}_{t_0} R_{\ell+k}(t)$ at each point $t_0\in\Delta$ there holds an estimate
$$
-{\rm ord}_{t_0}R_{\ell+k}(t)\leqslant k\,C,
$$
where $C>0$ is a constant that depends on the equation only.
\end{lemma}

{\bf Proof.} First of all let us observe the following differentiation rule:
$$
\delta: \; R_{\ell+k}(\ln x)\,x^k \;\mapsto\; x^k\,\Bigl(k+\frac{d}{dt}\Bigr)R_{\ell+k}(t)|_{t=\ln x},
$$
which implies
\begin{eqnarray*}
(\delta+\ell)^j: & & R_{\ell+k}(\ln x)\,x^k \;\mapsto\; x^k\,\Bigl(\ell+k+\frac{d}{dt}\Bigr)^jR_{\ell+k}(t)|_{t=\ln x}, \quad j=0,1,\ldots,n, \\
L(\delta):  & & R_{\ell+k}(\ln x)\,x^k \;\mapsto\; x^k\,\sum_{j=0}^na_j(t)\Bigl(\ell+k+\frac{d}{dt}\Bigr)^j
R_{\ell+k}(t)|_{t=\ln x}.
\end{eqnarray*}
Therefore substituting $\psi=\sum_{k=1}^{\infty}R_{\ell+k}(\ln x)x^k$ into the equation
$$
L(\delta)u=x\,M(x,\ln x, u,\delta u,\dots,\delta^n u)
$$
one obtains the equality the both sides of which are power-log series. Comparing the rational functions of $t=\ln x$, the coefficients at the first power of $x$ in these two series, one has
$$
\sum_{j=0}^na_j(t)\Bigl(\ell+1+\frac{d}{dt}\Bigr)^jR_{\ell+1}(t)=M(0,t,0,\dots,0)=b_{{\bf 0},0}(t).
$$

Let $\Delta$ be a set of poles of the finite collection of rational functions $a_0(t),\ldots,a_{n-1}(t)$,
$\{b_{{\bf q},\mu}(t)\}$ (not taking in account infinity). Then from general theory of linear differential equations it follows that all the finite poles of $R_{\ell+1}(t)$ belong to $\Delta$. For each point $t_0\in\Delta$ we define
$$
\gamma_{t_0}=\max_{j=0,1,\ldots,n}\bigl(j-{\rm ord}_{t_0}a_j(t)\bigr)\geqslant n.
$$
If this maximum is achieved on the coefficients $a_{j_1},\ldots,a_{j_p}$ then each of these coefficients has the Laurent expansion of the form
$$
a_{j_i}(t)=\frac{A_{j_i}}{(t-t_0)^{\gamma_{t_0}-j_i}}+\ldots, \qquad A_{j_i}\ne0, \quad i=1,\ldots,p,
$$
near the point $t=t_0$. Therefore, if
$$
R_{\ell+1}(t)=\frac{\alpha_1}{(t-t_0)^{\nu_1}}+\ldots, \qquad \alpha_1\ne0,
$$
in a punctured neighbourhood of $t_0$, then one has
$$
b_{{\bf 0},0}(t)=\sum_{j=0}^na_j(t)\Bigl(\ell+1+\frac{d}{dt}\Bigr)^jR_{\ell+1}(t)=\frac{\alpha_1P_{t_0}(-\nu_1)}
{(t-t_0)^{\nu_1+\gamma_{t_0}}}+\ldots
$$
in the same neighbourhood, where
$$
P_{t_0}(\lambda)=A_{j_p}\lambda(\lambda-1)\ldots(\lambda-j_p+1)+\ldots+
A_{j_1}\lambda(\lambda-1)\ldots(\lambda-j_1+1).
$$
Thus we have the following dichotomy for the order ${\rm ord}_{t_0}R_{\ell+1}(t)=-\nu_1$:
\begin{itemize}
\item either $P_{t_0}(-\nu_1)=0$, and then $\nu_1$ does not exceed the modules of the integer roots of
$P_{t_0}(\lambda)$;
\item or $P_{t_0}(-\nu_1)\ne0$, and then $\nu_1=-{\rm ord}_{t_0}b_{{\bf 0},0}(t)-\gamma_{t_0}<
-{\rm ord}_{t_0}b_{{\bf 0},0}(t)$.
\end{itemize}
Denoting by $C_1$ the maximum of the modules of the integer roots of $P_{t_0}(\lambda)$, and by $C_2$ the number
$$
C_2=\max_{{\bf q},\mu}\bigl(-{\rm ord}_{t_0}b_{{\bf q},\mu}(t)+q_1+2q_2+\ldots+nq_n\bigr),
$$
where the maximum is taken over a finite set of elements {\bf q}, $\mu$ involved in (\ref{sumM}), we put  $C=\max(C_1,C_2)$. Further we prove the assertion of the lemma  by induction with respect to $k$, which has been already checked for $k=1$.

Let us denote by $R_k^j(t)$ the rational function $\bigl(\ell+k+\frac{d}{dt}\bigr)^jR_{\ell+k}(t)$, $j=0,1,\ldots,n$ (in particular, $R_k^0=R_{\ell+k}$). Then
$$
(\delta+\ell)^j\psi=\sum_{k=1}^{\infty}R_k^j(\ln x)\,x^k,
$$
and consistently comparing the coefficients at each power $x^k$ (starting with $k=2$) in the equality for the two formal power-log series,
$$
L(\delta)\psi=x\,M(x, \ln x, \psi,\delta\psi,\dots,\delta^n\psi),
$$
we obtain nonhomogeneous linear ODEs satisfied by the corresponding $R_{\ell+k}(t)$, $k\geqslant2$:
\begin{equation}\label{lODE}
\sum_{j=0}^na_j(t)\Bigl(\ell+k+\frac{d}{dt}\Bigr)^jR_{\ell+k}(t)=\widetilde R_k(t),
\end{equation}
where $\widetilde R_k$ is a finite sum of the functions of the form
$$
b_{{\bf q},\mu}(t)\bigl(R^0_{k_1}\ldots R^0_{k_{q_0}}\bigr)\bigl(R^1_{l_1}\ldots R^1_{l_{q_1}}\bigr) \ldots\bigl(R^n_{m_1}\ldots R^n_{m_{q_n}}\bigr),
$$
$$
\sum_{i=1}^{q_0}k_i+\sum_{i=1}^{q_1}l_i+\ldots+\sum_{i=1}^{q_n}m_i\leqslant k-1,
$$
which necessarily has finite poles at points of the set $\Delta$ only, due to the inductive assumption. Therefore the function $R_{\ell+k}$ being a solution of the linear ODE (\ref{lODE}) may have finite poles also at points of this set only. Due to the same inductive assumption one additionally has
\begin{eqnarray*}
-{\rm ord}_{t_0} \bigl(R^0_{k_1}\ldots R^0_{k_{q_0}}\bigr) & \leqslant & (k_1+\ldots+k_{q_0})C, \\
-{\rm ord}_{t_0} \bigl(R^1_{l_1}\ldots R^1_{l_{q_1}}\bigr) & \leqslant & (l_1+\ldots+l_{q_1})C+q_1, \\
\ldots & \ldots & \ldots \\
-{\rm ord}_{t_0} \bigl(R^n_{m_1}\ldots R^n_{m_{q_n}}\bigr) & \leqslant & (m_1+\ldots+m_{q_n})C+nq_n,
\end{eqnarray*}
whence
$$
-{\rm ord}_{t_0}\widetilde R_k(t)\leqslant C_2+(k-1)C\leqslant k\,C.
$$

Finally, for the order ${\rm ord}_{t_0}R_{\ell+k}(t)=-\nu_k$ one has a similar dichotomy:
\begin{itemize}
\item either $P_{t_0}(-\nu_k)=0$, and then $\nu_k\leqslant C_1$;
\item or $P_{t_0}(-\nu_k)\ne0$, and then $\nu_k=-{\rm ord}_{t_0}\widetilde R_k(t)-\gamma_{t_0}<
-{\rm ord}_{t_0}\widetilde R_k(t)\leqslant k\,C$.
\end{itemize}
The lemma is proved.

\begin{lemma}\label{l3}
For $\ell>m$ large enough, the formal series $\psi$ is a unique power-log series solution $($with coefficients rational in $\ln x)$ of $(\ref{rODE})$. Additionally one has the following estimate for the order ${\rm ord}_{\infty}R_{\ell+k}(t)$:
$$
-{\rm ord}_{\infty}R_{\ell+k}(t)\leqslant k\,C,
$$
where $C>0$ is a constant depending on the equation only.
\end{lemma}

{\bf Proof.} Due to Remark 1, $a_n(t)\equiv1$ and therefore all the coefficients $a_j$ of the operator $L(\delta)$ can be presented in the following form in a neighbourhood of the point $t=\infty$:
$$
a_j(t)=a_{j,p}\,t^p+a_{j,p-1}\,t^{p-1}+\ldots\, \qquad j=0,1,\ldots,n,
$$
where $p\geqslant0$ is common for all $j$ and such that among the leading coefficients  $a_{j,p}$, $j=0,1,\ldots,n$, there is a nonzero one. Let
$$
R_{\ell+k}(t)=\beta_k\,t^{\mu_k}+\beta_k^1\,t^{\mu_k-1}+\ldots, \qquad \beta_k\ne0,
$$
in this neighbourhood, then
$$
a_j(t)\Bigl(\ell+k+\frac{d}{dt}\Bigr)^jR_{\ell+k}(t)=\beta_k\,a_{j,p}(\ell+k)^jt^{\mu_k+p}+\ldots,
$$
and hence for $\widetilde R_k$ from (\ref{lODE}) one has:
$$
\widetilde R_k(t)=\sum_{j=0}^na_j(t)\Bigl(\ell+k+\frac{d}{dt}\Bigr)^jR_{\ell+k}(t)=\beta_k\,P_{\infty}(\ell+k)t^{\mu_k+p}
+\ldots,
$$
where
$$
P_{\infty}(\lambda)=a_{n,p}\,\lambda^n+a_{n-1,p}\,\lambda^{n-1}+\ldots+a_{1,p}\,\lambda+a_{0,p}\not\equiv0.
$$
Choose $\ell>m$ in such a way that the polynomial
$$
P_{\infty}(\lambda)=a_{n,p}\lambda^n+a_{n-1,p}\lambda^{n-1}+\ldots+a_{1,p}\lambda+a_{0,p}
$$
does not have integer roots greater than $\ell$. Then
\begin{equation}\label{orderinf}
-{\rm ord}_{\infty}R_{\ell+k}(t)=\mu_k=-p-{\rm ord}_{\infty}\widetilde R_k(t)\leqslant-{\rm ord}_{\infty}\widetilde R_k(t)
\leqslant k\,C,
\end{equation}
if a constant $C$ is given by
$$
C=\max_{{\bf q},\mu}\bigl(-{\rm ord}_{\infty}b_{{\bf q},\mu}(t)\bigr).
$$
The last inequality is proved by induction similarly to what we did in Lemma \ref{l2} for the order of a finite pole, the only difference is that the differentiation increases the order of a finite pole by 1, while it lowers by 1 the order of pole at infinity.

If there existed another formal power-log series solution of (\ref{rODE}) besides $\psi$, then some of the equations (\ref{lODE}) would possess a rational solution different from $R_{\ell+k}$. In this case the corresponding homogeneous linear ODE
$$
\sum_{j=0}^na_j(t)\Bigl(\ell+k+\frac{d}{dt}\Bigr)^ju=0
$$
would possess a {\it nonzero} rational solution $R(t)$, which is impossible. Indeed, for $R(t)$ there holds an equality coming from (\ref{orderinf}),
$$
-p-{\rm ord}_{\infty}\sum_{j=0}^na_j(t)\Bigl(\ell+k+\frac{d}{dt}\Bigr)^jR(t)=-{\rm ord}_{\infty}R(t),
$$
whereas the zero function has an infinite order. The lemma is proved.
\medskip

{\bf Remark 2.}  Due to Lemmas \ref{l2}, \ref{l3} the formal series $\psi$ satisfying the equation (\ref{rODE}) can be presented in the form
\begin{equation}\label{psiseries}
\psi=\sum\limits_{k=1}^{\infty}\frac{P_k(\ln x)}{Q^k(\ln x)}\,x^k,
\end{equation}
where $Q(t)=(t-t_1)^C\ldots(t-t_s)^C$ and $P_k$'s are polynomials of degree $d_k\leqslant c_1k$ (a constant $C\in{\mathbb Z}_+$ can be chosen common for the both lemmas, then $c_1=(s+1)C$).

\section{Main estimates}\label{s3}

In this section we obtain main estimates for the norm of the coefficients of the series \eqref{psiseries} and its derivatives $\delta^j\psi$, which will be further used in the proof of Theorem \ref{th1}.

For simplicity let us redefine the operator $L(\delta)$ in the left hand side of (\ref{rODE}) as
\begin{equation}\label{Ld}
L(\delta)=\sum\limits_{j=0}^n a_j(\ln x)\,\delta^j,
\end{equation}
so that the coefficients $a_j(t)$ preserve their properties at the point $t=\infty$.

The operator $\delta$ transforms any function $x^k\,P(\ln x)/Q^{m}(\ln x)$, where $P$ is a polynomial, to the function
$$
x^k\,\frac{\bigl[kQ(t)-mQ'(t)+Q(t)(d/dt)\bigr]P(t)|_{t=\ln x}}{Q^{m+1}(\ln x)}=:x^k\,\frac{D_{k,m}(P)}{Q^{m+1}(\ln x)},
$$
therefore,
\begin{eqnarray*}
\delta^j: & & \frac{P_k(\ln x)}{Q^k(\ln x)}\,x^k\;\mapsto\;\frac{D_{k,k+j-1}\circ\ldots\circ D_{k,k+1}\circ D_{k,k}(P_k)}{Q^{k+j}(\ln x)}\,x^k, \quad j=1,\ldots,n,\\
L(\delta): & & \frac{P_k(\ln x)}{Q^k(\ln x)}\,x^k\;\mapsto\;x^k\,\sum_{j=0}^n\frac{a_j(\ln x)}{Q^{k+j}(\ln x)}\,
D_{k,k+j-1}\circ\ldots\circ D_{k,k+1}\circ D_{k,k}(P_k).
\end{eqnarray*}

For any function $f$ that is meromorphic at infinity and represented there by a Laurent series with a finite principal part,
$$
f(t)=\sum_{i\geqslant p}f_i\,t^{-i},
$$
let us define its finite norm as
\begin{equation}\label{norm}
\|f\|=\|f\|_r=\sum_{i\geqslant p}|f_i|\,r^{-i},
\end{equation}
where a sufficiently large $r$ will be defined later. Then for an arbitrary polynomial $P$ and its derivative $P'$ one has
$$
\|P'\|\leqslant (d/r)\|P\|, \qquad d=\deg P,
$$
therefore from standard properties of the norm and from the fact that in our case the norm of the product of two functions does not exceed the product of their norms, it follows that
\begin{equation}\label{eq4}
\|D_{k,m}(P)\|\leqslant (k+m+d)c_2\,\|P\|, \qquad c_2=c_2(r)>0.
\end{equation}
On other hand, from the relation
$$
(kQ-mQ')P=-QP'+D_{k,m}(P)
$$
one obtains
$$
P=\frac1{1-(m/k)(Q'/Q)}\Bigl(-\frac1k\,P'+\frac1{kQ}\,D_{k,m}(P)\Bigr).
$$
The function $Q'/Q$ is holomorphic at infinity and vanishes there, hence the function
$\bigl(1-(m/k)(Q'/Q)\bigr)^{-1}$ is holomorphic at infinity and takes the value 1 there.
Thus if the ratio $m/k$ is bounded (further we will be interested in the case $k\leqslant m\leqslant k+n-1$), choosing $r$ in an appropriate way one can achieve an estimate
\begin{equation}\label{norm1}
\Bigl\|\frac1{1-(m/k)(Q'/Q)}\Bigr\|\leqslant 2,
\end{equation}
uniform with respect to $k$ and $m$. Therefore,
$$
\|P\|\leqslant\frac{2d}{kr}\,\|P\|+\frac{c_3}{2k}\,\|D_{k,m}(P)\|, \qquad c_3=c_3(r)>0,
$$
whence it follows that the norm of any polynomial $P$ whose degree $d$ does not exceed the value $kr/4$
satisfies the estimate
\begin{equation}\label{eq5}
  \|P\|\leqslant\frac{c_3}k\,\|D_{k,m}(P)\|.
\end{equation}

{\it Choice of the number $r$ determining the norm $\|\cdot\|_r$.} Let us choose the value $r$ in such a way that
\begin{itemize}
\item for every $k\geqslant1$, $k\leqslant m\leqslant k+n-1$, the inequality (\ref{norm1}) holds;
\item for every $k\geqslant1$ the inequality $c_1k+(n-1)\deg Q\leqslant kr/4$ holds;
\item the inequalities $\|1/Q\|<1$, $\|Q^{n+1}\|>1$ hold.
\end{itemize}
(Note that the constant $c_1$, which has appeared in Remark 2, does not depend on $r$.) Further by
$\|\cdot\|=\|\cdot\|_r$ we will mean the norm (\ref{norm}) defined by this very number $r$.
\medskip

For $j=0,1,\ldots,n$ denote by ${\cal D}_{k,j}$ the differential operator
$$
{\cal D}_{k,j}=D_{k,k+j-1}\circ\ldots\circ D_{k,k+1}\circ D_{k,k}
$$
of order $j$  (with ${\cal D}_{k,0}={\rm id}$) and define
\begin{equation}\label{eq8}
{\cal L}_k=\sum_{j=0}^na_j(t)Q^{n-j}(t)\,{\cal D}_{k,j}.
\end{equation}
In the next two lemmas we obtain the main estimates which relate the norms of the polynomials ${\cal D}_{k,j}(P_k)$ with the norms of the polynomials $P_k$ and ${\cal D}_{k,n}(P_k)$, as well as the norm of ${\cal D}_{k,n}(P_k)$ with the norm of ${\cal L}_k(P_k)$.

\begin{lemma}\label{l4} For the norm of the polynomials ${\cal D}_{k,j}(P_k)$, $j=0,1,\ldots,n$, the following estimates hold:
\begin{equation}\label{eq6}
\|{\cal D}_{k,j}(P_k)\|\leqslant \tilde c_2\,k^j\,\|P_k\|,
\end{equation}
\begin{equation}\label{eq7}
\|{\cal D}_{k,j}(P_k)\|\leqslant \frac{\tilde c_3}{k^{n-j}}\,\|{\cal D}_{k,n}(P_k)\|,
\end{equation}
with some $\tilde c_2, \tilde c_3>0$ not depending on $k$.
\end{lemma}

{\bf Proof.} Recall that $\deg P_k=d_k\leqslant c_1k$. From the definition of the operator $D_{k,m}$ it follows that
$\deg D_{k,m}(P)=\deg P+\deg Q$, for an arbitrary polynomial $P$. Therefore
$\deg {\cal D}_{k,j-1}(P_k)=d_k+(j-1)\deg Q$, and one comes to the estimate (\ref{eq6}) for $j\geqslant1$ by a
$j$-multiple application of the inequality (\ref{eq4}), starting with $P={\cal D}_{k,j-1}(P_k)$ and $m=k+j-1$.

Similarly, $\deg {\cal D}_{k,j}(P_k)=d_k+j\deg Q\leqslant kr/4$ due to the second condition determining the choice of the value $r$, and the estimate (\ref{eq7}) for $j\leqslant n-1$ is obtained by a $(n-j)$-multiple application of the inequality (\ref{eq5}), starting with $P={\cal D}_{k,j}(P_k)$ and $m=k+j$. The lemma is proved.

\begin{lemma}\label{l5} For the norm of ${\cal D}_{k,n}(P_k)$ the following estimate holds:
\begin{equation}\label{eq9}
\|{\cal D}_{k,n}(P_k)\|\leqslant A\,\|{\cal L}_k(P_k)\|, \qquad k\geqslant 1,
\end{equation}
with some $A>0$ not depending on $k$.
\end{lemma}

{\bf Proof.} Since
$$
{\cal L}_k(P_k)={\cal D}_{k,n}(P_k)+a_{n-1}(t)Q(t){\cal D}_{k,n-1}(P_k)+\ldots+
$$
$$
+a_1(t)Q^{n-1}(t){\cal D}_{k,1}(P_k)+a_0(t)Q^n(t){\cal D}_{k,0}(P_k)
$$
(recall that $a_n=1$ due to Remark 1), from (\ref{eq7}) it follows that
$$
\|{\cal D}_{k,n}(P_k)\|\leqslant\|{\cal L}_k(P_k)\|+\frac{\alpha_1}k\,\|{\cal D}_{k,n}(P_k)\|+\ldots+
\frac{\alpha_n}{k^n}\,\|{\cal D}_{k,n}(P_k)\|,
$$
with some $\alpha_1,\ldots,\alpha_n>0$ not depending on $k$. Therefore
$$
\Bigl(1-\frac{\alpha_1}k-\ldots-\frac{\alpha_n}{k^n}\Bigr)\|{\cal D}_{k,n}(P_k)\|\leqslant\|{\cal L}_k(P_k)\|.
$$
The value in brackets is positive and separated from zero beginning with some $k$, which implies the lemma.

\section{Constructing an equation with a majorant solution}\label{s4}

In this section using the reduced ODE \eqref{rODE} with the formal solution \eqref{psiseries} we construct a functional equation
\begin{equation}\label{eq10}
\sigma\, U=x\,\widetilde{M}(x,\ln x,U), \qquad \sigma=1/(\tilde{c}_3 A),
\end{equation}
$\widetilde M$ being a polynomial of its variables, with a formal solution $U=\tilde{\psi}$ in the form of a {\it Dulac series}
\begin{equation}\label{eq11}
\tilde{\psi}=\sum_{k=1}^{\infty}\widetilde{P}_k(\ln x)\,x^k,\qquad \widetilde{P}_k\in\mathbb{R}_+[t],
\end{equation}
and prove that the series \eqref{eq11} is \textit{majorant} for the series \eqref{psiseries}, that is, for the norms of the coefficients of these series the following inequalities will hold:
\begin{equation}\label{eq12}
\displaystyle \left\|\frac{P_k}{Q^k}\right\|\leqslant \tilde c^k\,\|\widetilde{P}_k\|, \qquad k\geqslant 1,
\end{equation}
where $\tilde c=\|Q^{1+n}\|>1$.

A polynomial $\widetilde{M}$ is defined in the following way. After introducing the polynomial $Q$, the expression
(\ref{sumM}) determining $M(x,\ln x,u,\delta u,\ldots,\delta^nu)$ can be presented as a finite sum of monomials of the form
\begin{equation}\label{Msummand}
\alpha\,x^{\mu}\,\frac{(\ln x)^{\nu}}{Q(\ln x)}\, u^{q_0}(\delta u)^{q_1}\ldots(\delta^nu)^{q_n}, \qquad
\alpha=\alpha_{{\bf q},\mu,\nu}\in\mathbb C.
\end{equation}
To define the polynomial $\widetilde{M}(x,\ln x, U)$ we change each such a monomial in $M$ to the corresponding monomial
\begin{equation}\label{tMsummand}
|\alpha|\,x^{\mu}(\ln x)^{\nu}\,U^{q_0}U^{q_1}\ldots U^{q_n}.
\end{equation}

Let us prove that the equation \eqref{eq10} obtained {\it via} such a construction, is majorant for the ODE \eqref{rODE} in the following sense.

\begin{lemma}\label{l6}
{\rm (i)} The equation \eqref{eq10} possesses a unique formal solution in the form of the Dulac series \eqref{eq11}, furthermore $\widetilde{P}_k\in\mathbb{R}_+[t]$ are polynomials with non-negative real coefficients;

{\rm (ii)} The series \eqref{eq11} is majorant for the formal series \eqref{psiseries} satisfying the equation \eqref{rODE}, that is, for the norms of the corresponding coefficients $\widetilde{P}_k$ and $P_k/Q^k$ of these series the inequality  \eqref{eq12} holds.
\end{lemma}

{\bf Proof.} (i) As follows from \eqref{eq10}, every coefficient $\widetilde P_k$ is uniquely determined {\it via} the previous ones $\widetilde P_1,\ldots,\widetilde P_{k-1}$ by a polynomial expression, where
$\widetilde P_1(t)=(1/\sigma)\widetilde M(0,t,0)\in{\mathbb R}_+[t]$. The coefficients of such a polynomial expression are positive, as follows from the form \eqref{tMsummand} of the summands of the polynomial $\widetilde M(x,t,U)$. This proves the part (i).
\medskip

{\bf Remark 3.} For any polynomials $P,Q\in{\mathbb R}_+[t]$ one has $\|P+Q\|=\|P\|+\|Q\|$ and $\|PQ\|=\|P\|\cdot\|Q\|$.
\medskip

(ii) Since the action of the linear operator (\ref{Ld}) to each term of the series \eqref{psiseries} has the form
$$
L(\delta):  \frac{P_k(\ln x)}{Q^k(\ln x)}\,x^k\;\mapsto\;x^k\,\sum_{j=0}^n\frac{a_j(\ln x)}{Q^{k+j}(\ln x)}\,
{\cal D}_{k,j}(P_k)=x^k\,\frac{{\cal L}_k(P_k)}{Q^{k+n}(\ln x)},
$$
the polynomials $P_k$ satisfy nonhomogeneous linear ODEs,
\begin{equation}\label{ODEforP}
\frac{\mathcal{L}_k(P_k(t))}{Q^{k+n}(t)}=B_k(t), \qquad k=1,2,\ldots,
\end{equation}
where $B_1(t)=M(0,t,0,\ldots,0)$ and the other $B_k$'s are also rational functions, which will be explicitly described below.

Let us first prove the estimates $\|{\cal L}_k(P_k)\|\leqslant \sigma\,\tilde c^k\|\widetilde P_k\|$, $k\geqslant1$, using the induction with respect to $k$. For $k=1$ one has:
$$
\|B_1\|=\|M(0,t,0,\ldots,0)\|\leqslant\|\widetilde{M}(0,t,0)\|=\sigma\|\widetilde P_1\|
$$
(the inequality $\|M(0,t,0,\ldots,0)\|\leqslant\|\widetilde{M}(0,t,0)\|$ is provided by construction and by the estimate $\|1/Q\|<1$), whence
$$
\|{\cal L}_1(P_1)\|=\|Q^{1+n}\,B_1\|\leqslant \sigma\,\tilde c\|\widetilde P_1\|.
$$

To obtain analogous estimates for all $k\geqslant2$ let us describe in detail the function $B_k$ appearing in (\ref{ODEforP}) and the expression for the coefficient $\widetilde P_k$ mentioned in the proof of the part (i). Put $P_k^j=\mathcal{D}_{k,j}(P_k)$ (in particular, $P_k^0=P_k$).

Since
$$
\delta^j\psi=\sum_{k=1}^{\infty}\frac{P_k^j(\ln x)}{Q^{k+j}(\ln x)}\,x^k, \qquad
\tilde\psi=\sum_{k=1}^{\infty}\widetilde P_k(\ln x)\,x^k,
$$
coming back to the expression \eqref{Msummand} for $u=\psi$ one concludes that $B_k(t)$ is a sum of functions of the form
$$
Q^{-1}(t)\cdot \alpha\,t^{\nu}\,\Bigl(\frac{P^0_{k_1}}{Q^{k_1}}\ldots\frac{P^0_{k_{q_0}}}{Q^{k_{q_0}}}\Bigr) \Bigl(\frac{P^1_{l_1}}{Q^{l_1+1}}\ldots\frac{ P^1_{l_{q_1}}}{Q^{l_{q_1}+1}}\Bigr)\ldots
\Bigl(\frac{P^n_{m_1}}{Q^{m_1+n}}\ldots\frac{P^n_{m_{q_n}}}{Q^{m_{q_n}+n}}\Bigr),
$$
where $\sum_{i=1}^{q_0}k_i+\sum_{i=1}^{q_1}l_i+\ldots+\sum_{i=1}^{q_n}m_i=k-1-\mu$, $\mu$ being the degree of $x$ in the corresponding monomial \eqref{Msummand}. Therefore the function ${\cal L}_k(P_k)=Q^{k+n}\,B_k$ is a sum of functions of the form
\begin{equation}\label{Qsummand}
Q^{n+\mu-q_1-2q_2-\ldots-nq_n}(t)\cdot \alpha\,t^{\nu}\,(P^0_{k_1}\ldots P^0_{k_{q_0}})(P^1_{l_1}\ldots P^1_{l_{q_1}})\ldots(P^n_{m_1}\ldots P^n_{m_{q_n}}),
\end{equation}
whereas due to \eqref{tMsummand} with $U=\tilde\psi$, $\sigma\widetilde P_k(t)$ is a sum of the corresponding polynomials
\begin{equation}\label{tQsummand}
|\alpha|\,t^{\nu}\,(\widetilde P_{k_1}\ldots\widetilde P_{k_{q_0}})
(\widetilde P_{l_1}\ldots\widetilde P_{l_{q_1}})\ldots(\widetilde P_{m_1}\ldots\widetilde P_{m_{q_n}}).
\end{equation}

The norm of the function \eqref{Qsummand} does not exceed the product
$$
\tilde c^{\mu+1}\,|\alpha|\cdot r^\nu \cdot \|P^0_{k_1}\|\ldots\|P^0_{k_{q_0}}\|\cdot
\|P^1_{l_1}\|\ldots\|P^1_{l_{q_1}}\|\cdot\ldots\cdot\|P^n_{m_1}\|\ldots\|P^n_{m_{q_n}}\|.
$$
By the inductive assumption and Lemmas \ref{l4}, \ref{l5}, with the use of the relation $\tilde c_3A=1/\sigma$, we can estimate each factor $\|P_s^j\|$, $s<k$, in the last product:
$$
\|P_s^j\|=\|{\cal D}_{s,j}(P_s)\|\leqslant{\tilde{c}_3}\,\|{\cal D}_{s,n}(P_s)\|\leqslant{\tilde{c}_3}A\,
\|\mathcal{L}_s(P_s)\|\leqslant {\tilde c}^s\|\widetilde{P}_s\|.
$$
It follows that the norm of (\ref{Qsummand}) does not exceed the norm of (\ref{tQsummand}) multiplied by $\tilde c^k$
(recall that the norm of the product of polynomials with non-negative real coefficients equals the product of the norms of these polynomials), thus we come to the required auxiliary estimate $\|{\cal L}_k(P_k)\|\leqslant\sigma\,\tilde c^k\,
\|\widetilde P_k\|$ (again, the norm of the sum of polynomials with non-negative real coefficients equals the sum of their norms). Finally the part (ii) follows from the chain of inequalities
\begin{eqnarray*}
\left\|\frac{P_k}{Q^k}\right\|\leqslant \|P_k\|\leqslant\tilde c_3\, \|{\cal D}_{k,n}(P_k)\|\leqslant  {\tilde{c}_3}A\,
\|\mathcal{L}_k(P_k)\|\leqslant \tilde c^k\,\|\tilde{P}_k\|.
\end{eqnarray*}
The lemma is proved.

\section{Finishing the proof of Theorem \ref{th1}}\label{s5}

As follows from \cite[\S5]{GG}, the series $\sum_{k\geqslant1}\widetilde P_k(|\ln x|)\,x^k$ converges absolutely in any open sector $S$ of sufficiently small radius with vertex at the origin and of opening less than $2\pi$.

Recall that $\deg P_k=d_k\leqslant c_1k$, therefore for sufficiently small $|x|$ one has
$$
\left|\frac{P_k(\ln x)}{Q^k(\ln x)}x^k\right|\leqslant |P_k(\ln x)|\cdot |x|^k\leqslant \|P_k\|\cdot|\ln x|^{d_k}\,|x|^k\leqslant \|P_k\|\cdot|x|^{\varepsilon k}\leqslant \left\|\frac{P_k}{Q^k}\right\|\cdot\|Q\|^k\,|x|^{\varepsilon k},
$$
for some $0<\varepsilon<1$. Denoting $\alpha=\tilde c\,\|Q\|>0$ and taking into account the estimate \eqref{eq12},
for sufficiently small $|x|$ we obtain
$$
\left|\frac{P_k(\ln x)}{Q^k(\ln x)}x^k\right|\leqslant\|\widetilde P_k\|\cdot\alpha^k\,|x|^{\varepsilon k}\leqslant
\widetilde P_k\bigl(|\ln(\alpha x^{\varepsilon})|\bigr)\,|\alpha x^{\varepsilon}|^k,
$$
whence it follows the uniform convergence of the series \eqref{psiseries}, together with that of (\ref{series}),
in any open sector $S$ of sufficiently small radius with vertex at the origin and of opening less than $2\pi$. This finishes the proof of Theorem \ref{th1}.
\medskip

Further we propose several examples of formal power-log series solutions of algebraic ODEs and discuss possible applications of Theorem 1 to the study of their convergence.
\medskip

{\bf Example 1.} Consider the sixth Painlev\'e equation
$$
y''=\frac{(y')^2}{2}\left(\frac{1}{y}+\frac{1}{y-1}+\frac{1}{y-x}\right)-y'\left(\frac{1}{x}
+\frac{1}{x-1}+\frac{1}{y-x}\right)+
$$
\begin{equation}
+\frac{y(y-1)(y-x)}{x^2(x-1)^2}\left[a+b\frac{x}{y^2}+c\frac{x-1}{(y-1)^2}+
d\frac{x(x-1)}{(y-x)^2}\right],\label{p6}
\end{equation}
where  $a,$ $b,$ $c,$ $d$ are complex parameters. If $1-2d+2b\neq 0$ then the equation \eqref{p6} possesses  a one-parameter family of formal Dulac series solutions \cite{BG}, \cite{Guz}, \cite{Sh},
\begin{equation}\label{p62}
\varphi=\sum\limits_{k=1}^{\infty}P_k(\ln x)x^k, \qquad P_k\in{\mathbb C}[t],
\end{equation}
where $P_1(\ln x)=\frac14(1-2d+2b)(\ln x+C)^2+2b/(1-2d+2b)$, the other polynomials $P_k$ are uniquely determined and the degree of each $P_k$ does not exceed $2k$. More precisely, the polynomiality of $P_k$'s is established by S.\,Shimomura \cite{Sh}, who has also proved the convergence of the series (\ref{p62}) for small $|x|$.

One can also establish convergence with the use of Theorem \ref{th1}, rewriting the equation (\ref{p6}) in a polynomial form $F(x,y,\delta y,\delta^2 y)=0$ and directly checking that the series $\partial F/\partial y_j(x,\Phi)$
begin with $({\rm polynomial\; in} \ln x)\cdot x^2$. Note that although the formal solution (\ref{p62}) has the form of a Dulac series, the theorem on convergence from the paper \cite{GG} formally cannot be applied here, since its assumptions require that the leading term of the partial derivative $\partial F/\partial y_2(x,\Phi)$ does not contain logarithm.
\medskip

{\bf Example 2.} In this example we describe a situation in which a formal power-log series solution of an algebraic ODE that satisfies the sufficient condition of convergence from Theorem \ref{th1} is produced from a formal {\it Dulac series} solution of another algebraic ODE that satisfies the sufficient condition of convergence from \cite{GG}, {\it via} a rational transformation of the unknown.

An equation
$$
\widetilde F(x,w,\delta w):=(\delta w)^2-4x^2w^3-1=0
$$
possesses a formal Dulac series solution
$$
\tilde\varphi=\ln x+\sum\limits_{k=1}^{\infty}P_{2k}(\ln x)x^{2k}, \qquad P_{2k}\in\mathbb{C}[t],
$$
where each $P_{2k}$ satisfies a linear ODE with constant coefficients and right hand side polynomially depending on the previous $P_2,\ldots,P_{2k-2}$:
$$
2\Bigl(\frac{d}{dt}+2k\Bigr)P_{2k}(t)=B_{2k}(t,P_2,\dots,P_{2k-2}), \qquad B_2(t)=4t^3
$$
(whence the polynomiality of $P_{2k}$'s follows). One has
\begin{eqnarray*}
\frac{\partial\widetilde F}{\partial w_1}(x,\tilde\varphi,\delta\tilde\varphi)&=&2\,\delta\tilde\varphi=2+\dots,\\
\frac{\partial\widetilde F}{\partial w_0}(x,\tilde\varphi,\delta\tilde\varphi)&=&-12x^2\tilde\varphi^2=(-12\ln^2x)x^2+\dots,
\end{eqnarray*}
hence $\tilde\varphi$ converges for sufficiently small $|x|$, due to \cite{GG}.

Thus an equation
$$
F(x,y,\delta y):=(\delta y)^2-4x^2y-y^4 =0
$$
obtained from the initial one by the transformation $y=1/w$ possesses a {\it convergent} power-log series solution
$$
\varphi=1/\tilde\varphi=\frac1{\ln x}+\sum\limits_{k=1}^{\infty}R_{2k}(\ln x)x^{2k}, \qquad R_{2k}\in\mathbb{C}(t).
$$
The assumption of Theorem \ref{th1} is expectingly fulfilled here, since
\begin{eqnarray*}
\frac{\partial F}{\partial y_1}(x,\varphi,\delta\varphi)&=&2\delta\varphi=-\frac{2}{\ln^2x}+\dots,\\
\frac{\partial F}{\partial y_0}(x,\varphi,\delta\varphi)&=&-4x^2-4\varphi^3=-\frac{4}{\ln^3x}+\dots.
\end{eqnarray*}
\medskip

Thus we see that Example 2 is the illustration of a simplest situation in which a formal power-log series solution $\varphi$ of (\ref{ODE}) that satisfies the sufficient condition of convergence from Theorem \ref{th1} (that is, when the linear operator $L(\delta)$ in the reduced equation (\ref{rODE}) obtained by $\varphi$ is of order $n$ exactly) is produced from a formal {\it Dulac series} solution $\tilde\varphi$ of another algebraic ODE {\it via} a rational transformation of the unknown, such that $\widetilde L(\delta)$ in the reduced ODE obtained by $\tilde\varphi$, is a linear diffrential operator of order $n$ with {\it constant} coefficients (thus the convergence of $\tilde\varphi$ follows from \cite{GG} in this case). In such a situation, the pairs $(\varphi,L)$ and $(\tilde\varphi,\widetilde L)$ are naturally considered equivalent and the convergence of $\varphi$ already follows from that of $\tilde\varphi$.

Leaving the case described in Example 2, where is almost no difference between power-log series and Dulac series solutions of an algebraic ODE, we note that there are a whole lot of less studied (but more general) examples, in which the first coefficients of a formal power-log transseries solution are rational functions of logarithm, and the rationality of the other coefficients having the form of formal Laurent series in $1/\ln x$ with a finite main part cannot be proved yet (though in most of such examples the computation of the first $N$ coefficients, with $N$ big enough, allows one to expect the rationality of the rest coefficients). Let us consider an example of such a general situation.
\medskip

{\bf Example 3.} An equation
$$
(\delta^2 y)^2-2(\delta y)^3-2x^2y-2x=0
$$
possesses a formal power-log transseries solution
$$
\varphi=-\frac2{\ln x}+\sum\limits_{k=1}^{\infty}\widehat R_k(\ln x)x^k, \qquad\widehat R_k\in\mathbb{C}((1/t)).
$$
The computation gives us the rationality of the first coefficients $\widehat R_k$: for example, $\widehat R_1\in(1/t^2){\mathbb C}[t]$, $\widehat R_2\in(1/t^3){\mathbb C}[t]$, $\widehat R_3\in(1/t^4){\mathbb C}[t]$. In the case of the rationality of all the $\widehat R_k$'s, the convergence of $\varphi$ would evidently follow from Theorem 1.
\medskip

{\bf Concluding remark.} Typical situations presented in the above two examples give rise to the following questions. {\it Does the rationality of the first coefficients of a formal power-log transseries solution of an algebraic ODE imply that the rest coefficients are also rational functions of logarithm? Among convergent power-log series solutions, are there those which cannot be obtained from Dulac series solutions {\rm via} a rational transformation of the unknown, and if yes, how can such cases be destinguished?}

\end{document}